\documentclass[3p,times]{article}
\usepackage{amsmath}
\usepackage{amsfonts}
\usepackage{mathrsfs}
\usepackage{amssymb}
\usepackage{times}
\usepackage{multirow}
\usepackage{geometry}
\usepackage{subcaption}
\usepackage{xcolor}
\usepackage{authblk}
\usepackage{graphicx}
\newtheorem{theorem}{Theorem}
\newtheorem{proposition}[theorem]{Proposition}
\newtheorem{remark}[theorem]{Remark}

\newtheorem{corollary}[theorem]{Corollary}
\newtheorem{lemma}[theorem]{Lemma}

\title{On some mixing properties of copula-based Markov chains}
\author[1]{Martial Longla\thanks{mlongla@olemiss.edu}}
\author[2]{Mous-Abou Hamadou.\thanks{mousabouhamadou@gmail.com}}
\author[3]{Seraphin Isidore Ngongo.\thanks{ngongoisidore@gmail.com}}
\affil[1]{University of Mississippi, Department of mathematics}
\affil[2]{University of Maroua, Department of mathematics}
\affil[3]{University of Yaounde I, ENS, Department of Mathematics}
\counterwithin{equation}{section}
\counterwithin{theorem}{subsection}
\begin{document}
\maketitle

\begin{abstract}
This paper brings some insights of $\psi'$-mixing, $\psi^*$-mixing and $\psi$-mixing for copula-based Markov chains and the perturbations of their copulas. We provide new tools to check Markov chains for $\psi$-mixing or $\psi'$-mixing, and also show that perturbations of $\psi'$-mixing copula-based Markov chains are $\psi'$-mixing while perturbations of $\psi$-mixing Markov chains are not necessarily $\psi$-mixing markov chains, even when the perturbed copula is $\psi-mixing$. Some examples of copula families are considered. A statistical study is provided to emphacize the impact of perturbations on copula-based Markov chains. Moreover, we provide a correction to a statement made in Longla and al. (2021) on $\psi$-mixing.
\end{abstract}

\textit{Key words}: Perturbations of copulas, mixtures of copulas, convex combinations of Copulas,
Mixing rates, Lower-psi mixing.

\textit{Mathematical Subject Classification} (2000): 62G08, 62M02, 60J35\bigskip

\section{introduction}
Modelling dependence among variables or factors in economics, finance, risk management and other applied fields has benefited over the last decades from the study of copulas. Copulas, these multivariate cumulative distributions with uniform marginals on $[0,1]^{n}$, have been widely used as strength of the dependence between variables. Sklar (1959) first showed that by rescalling the effect of marginal distributions, one obtains a copula from the joint distribution of random variables. This rescalling implies that when variables are transformed using increasing functions, the copula of their transformations remains same as that of the original variables. For many dependence coefficients, this copula is all that affects the computations (random vectors with common copulas have common dependence coefficients).  This justifies why dealing with the uniform distribution as stationary distribution of a Markov chain is same as studying a Markov chain with any absolutely continuous stationary distribution. Following the ideas of Durante and al. (2013), Longla and al. (2021) and Longla and al. (2022) have considered the perturbation method that adds to a copula an extra term called perturbation. They also considered other classes of modifications and their impact on the dependence structure as studied by Komornik and al. (2017). The long run impact of such perturbations on the dependence structure and the measures of association was investigated. In fact,They have investigated the impact of perturbations of copulas on the mixing structure of the Markov chains that they generate. The case was presented for $\rho$-mixing, $\alpha$-mixing, $\psi$-mixing and $\beta$-mixing in Longla and al. (2021) and Longla and al. (2022). 
\subsection{Facts about Copulas}
The definition of a 2-copula and related topics can be found in Nelsen (2006). 2-copulas are in  general just referred to as copulas when there is no reason for confusion. We will follow this assumption throughout this paper.  A function $C: [0,1]^{2}\rightarrow [0,1]$ is called a bivariate copula if it satisfies the following conditions:
\begin{enumerate}
	\item[i.] $C(0,x)=C(x,0)=0$ (meaning that $C$ is grounded); \item[ii.]$C(x,1)=C(1,x)=x, \forall x\in [0,1]$ (meaning each coordinate is uniform on [0,1]);
\item[iii.] $C(a,c)+C(b,d)-C(a,d)-C(b,c)\geq 0, \forall\
[a,b]\times[c,d]\subset [0,1]^{2}.	$
\end{enumerate}
The last condition basically states that the porbability of any rectangular subset of $[0,1]\times [0,1]$ is non-negative. This is an obvious condition, given that $C(x,y)$ is a cumulative probability distribution function on $[0,1]\times [0,1]$. The first condition states that the probability of any rectangle that doesn't cross $[0,1]\times [0,1]$ is equal to 0 (this covers that fact that such a rectangle doesn't intersect the support of the distribution function). The second condition basically asserts that the marginal distribution is uniform on $[0,1]$ for each of the components of the considered vector. 

Darsaw and al. (1992) derived the transition probabilities for stationary Markov chains with uniform marginals on $[0,1]$ as $P(X_{n}\in (-\infty,x]|X_{n-1}=x)=C_{,1}(x,y), \forall n\in\mathbb{N}$, where $C_{,i}(x,y)$ denotes the derivative of $C(x,y)$ with respect to the $i^{th}$ variable. This property has been used by many authors to establish mixing properties of copula-based Markov chains. We can cite Longla (2015), Longla (2014), Longla and Peligrad (2012) who provided some results for reversible Markov chains, Beare (2010) who presented results for $\rho$-mixing among others.  

It's been shown in the literature (see Darsow and al. (1992) and the references therein) that if $(X_1, \cdots, X_n)$ is a Markov chain with consecutive copulas $(C_1, \cdots, C_{n-1})$, then the fold product given by 
$$C(x,y)=C_1*C_2 (x, y)=\int^1_0 C_{1,2}(x, t)C_{2,1}(t, y)dt$$ is the copula of $(X_1,X_3)$ and the $\star$-product given by
$$ C(x,y,z)=C_1\star C_2 (x, y,z)=\int_0^y C_{1,2}(x, t)C_{2,1}(t, z)dt$$ is the copula of $(X_1,X_2,X_3)$. The $n$-fold product of $C(x,y)$ denoted $C^n(x,y)$ is defined by the recurrence $C^{1}(x,y)=C(x,y)$, $$C^{n}(x,y)=C^{n-1}*C(x,y).$$ 
The most popular copulas are $\Pi(u,v)=uv$  (the independent copula), the Hoeffding lower and upper bounds  $W(u,v)=\max(u+v-1,0)$ and $M(u,v)=\min(u,v)$ respectively.  Convex combinations of copulas $\{C_1(x,y), \cdots, C_k(x,y)\}$ defined by $\displaystyle \{ C(x,y)=\sum_{j=1}^{k}a_j C_j(x,y), 0\leq a_j, \sum_{j=1}^{k} a_j=1\}$ are also copulas. For any copula $C(x,y)$, there exists a unique representation $C(x, y) = AC(x, y) + SC(x, y)$, where
$AC(x, y)$ is the absolute continuous part of $C(x, y)$ and $SC(x, y)$ is the singular part of the copula $C(x,y)$. $AC(x,y)$ induces on $[0,1]^2$ a measure $P_c$ defined
on borel sets by
$$\displaystyle P_c(A\times B)=\int_A\int_B c(x,y)dxdy\quad
\text{and} \quad P(A\cap B)=P_c(A\times B)+SC(A\times B), \quad \text{(see Longla (2015)).}$$
An absolutely continuous copula is one that has singular part $SC(x,y)=0$ and a singular copula is one that has absolutely continuous part $AC(x,y)=0$. This work is concerned mostly by absolutely continuous copulas and mixing properties of the Markov chains they generate.

\subsection{Mixing coefficients of interest}
 The mixing coefficients of interest in this paper are $\psi'$ and $\psi$. The $\psi$-mixing condition has its origin in the paper by Blum and al. (1963). They studied a different condition (“$\psi$*-mixing”) similar to this mixing coefficient. They showed that for Markov chains satisfying their condition, the mixing rate is exponential. The coefficient took its present form in the paper of Philipp (1969). For examples of mixing sequences, see Kesten and O'Brien (1976), who showed that in general, the mixing rate could be arbitrarily slow, a large class of mixing rates can occur for stationary $\psi$-mixing. The general definitions of these mixing coefficients are as follows. Given any $\sigma$-fields $\mathscr{A}$ and $\mathscr{B}$ and a defined probability measure $P$,
$$\psi(\mathscr{A},\mathscr{B})=\sup_{B\in \mathscr{B}, A\in \mathscr{A}, P(A)\cdot P(B)>0 }\frac{|P(A\cap B)-P(A)P(B)|}{P(A)P(B)},$$
$$\psi'(\mathscr{A},\mathscr{B})=\inf_{B\in \mathscr{B}, A\in \mathscr{A}, P(A)>0}\frac{P(B\cap A)}{P(A)P(B)}, \quad \text{and}\quad \psi^*(\mathscr{A},\mathscr{B})=\sup_{B\in \mathscr{B}, A\in \mathscr{A}, P(A)\cdot P(B)>0 }\frac{P(A\cap B)}{P(A)P(B)}.$$ 
In case of stationary copula-based Markov chains generated by an absolutely continuous copula, the $\psi'$-mixing dependence coefficient takes the form
\begin{center}
$\psi'_n(C)=\underset{\underset{ \lambda(A)\lambda(B)>0}{ A,B\in \mathscr{B}}}{\inf}\dfrac{\int_A\int_B c_n(x,y)dxdy}{\lambda(A)\lambda(B)},$
\end{center} where  $c_n(x,y)$ is the density of the of $C^n(x,y) $ and $\lambda$ is the Lebesgue measure on $I=[0,1]$.  For every positive integer $n$, let $\mu_n$ be the measure induced by the distribution of $(X_0, X_n)$. Let $\mu$ be the measure induced by the stationary distribution of the Markov chain and $\mathscr{B}$ the $\sigma$-algebra generated by $X_0$. The $\psi'$-mixing dependence coefficient takes the form
\begin{center}
$\psi_n'(C)=\underset{A,B\in \mathscr{B}, \mu(A).\mu(B)>0~}{\inf}\dfrac{\mu_n(A\times B)}{\mu(A)\mu(B)}$, \quad and \quad $\psi_n^*(C)=\underset{A,B\in \mathscr{B}, \mu(A).\mu(B)>0~}{\sup}\dfrac{\mu_n(A\times B)}{\mu(A)\mu(B)}$
\end{center}

\subsection{About perturbations}
In applications, knowing approximately a copula $C(u,v)$ appropriate to the model of the observed data, minor perturbations of $C(u,v)$ are considered. Komornik and al. (2017) have investigated some perturbations that were introduced by Mesiar and al. (2015). These perturbations were also considered by Longla and al. (2021) and (2022). 
Perturbations that we consider in this work have been studied by many authors. Sheikhi et al. (2020) looked at the perturbations of copulas via modification of the random variables that the copulas are used to represent the dependence structure of. Namely, they perturbed the copula of $(X,Y)$ by looking at the copula of $(X+Z, Y+Z)$ for some $Z$ independent of $(X,Y)$ that can be considered as noise. Mesiar and al. (2019) worked on the perturbations induced by modification of one of the random variables of the pair. Namely, the copula of $(X,Y)$ was perturbed to obtain the copula of $(X+Z, Y)$. In this work, we look at he impact of perturbations on $\psi$-mixing and $\psi'$-mixing. We provide theoretical proofs and a simulation study that justifies the importance of the study of perturbations and their impact on estimation problems. This is done through the central limit theorem that varies from one kind of mixing structure to another and is severely impacted by perturbations in the case of $\psi$-mixing for instance.

\subsection{Structure of the paper}
This paper consists of six sections, each of which concern a specific topic of interest and are structured as follows. Introduction in Section 1 is divided into several parts. Facts about copulas are introduced in subsection 1.1, mixing coefficient of interest ($\psi'$-mixing and $\psi$-mixing)  are defined in subsection 1.2 and Subsection 1.3 is dedicated to facts about perturbation of copulas. Section 2 is devoted to the impact of perturbations on $\psi'$-mixing and $\psi$-mixing copula-based Markov chains, addressing $\psi'$-mixing in Subsection 2.1 and $\psi$-mixing in Subsection 2.2. We emphasize on the fact that perturbations of $\psi'$-mixing copula-based Markov chains are $\psi'$-mixing while perturbations of $\psi$-mixing Markov chains are not necessarily $\psi$-mixing, even when the perturbed copula is $\psi$-mixing. We present here the case of $\psi^*$-mixng. This section ends by an explicit example showing this fact. In Section 3 we provide some graphs to show the effect of perturbations. In Section 4, we showcase a simulation study to emphasize the importance of this topic. Comments on the paper's results and their relationship with current state of art are  presented in Section 5 and Section 6 provided the proofs of our main results. Throughout this work $\psi_n(C)$ is replaced by $\psi_n$ when there is no reason for confusion. 

\section{Facts about $\psi'$-mixing and $\psi$-mixing}
\subsection{All about $\psi'$-mixing}
A result of Bradley (1983) states the following
\begin{theorem}\label{Brad}
For any strictly stationary Markov chain, either $\psi'_n\to 1$ as $n\to\infty$ or $\psi'_n=0$ $\forall n\in\mathbb{N}$.
\end{theorem}
Based on this result, we show the following.
\begin{theorem}\label{TheoL} Let $\lambda$ be the Lebesgue measure on $[0,1]$.
If the copula $C(u,v)$ of the stationary Markov chain $(X_k, k\in\mathbb{N})$ is such that the density of its absolutely continuous part $c(u,v)\ge \varepsilon_1(u)+\varepsilon_2(v)$ on a set of Lebesgue measure $1$ and $\displaystyle \inf_{A\subset I}\frac{\int_{A}\varepsilon_1d\lambda}{\lambda(A)}>0$ or $\displaystyle \inf_{A\subset I}\frac{\int_{A}\varepsilon_2 d\lambda}{\lambda(A)}>0$, then the Markov chain is $\psi'$-mixing.
\end{theorem}
Theorem \ref{TheoL} is an extension of Theorem 2.5 of Longla (2014). It extends the result from $\rho$-mixing to $\psi'$-mixing. Longla and al. (2021) state that for a copula $C$ perturbed by means of the independence copula $\Pi$, the following result holds.
\begin{theorem} \label{copL} The perturbed copula with parameter $\theta$ has the following properties:
\begin{equation}
C_{\theta,\Pi}^{n}(u,v)=(1-\theta)^nC^{n}(u,v)+(1-(1-\theta)^n)uv.
\end{equation}
\end{theorem}
As a result of Theorem  \ref{copL}, following Longla (2015), based on the fact that the density of the copula  $C_{\theta,\Pi}^{n}(u,v)$ is bounded away from zero on a set of Lebesgue Measure $1$, we can conclude  the following:
\begin{corollary}
 $C_{\theta,\Pi}^{n}(u,v)$ generates lower $\psi$ mixing stationary Markov chains. 
\end{corollary}
In general, for any convex combination of copulas, the following result holds.
\begin{theorem}\label{TheoL2}
For any set of copulas $C_1(u,v)\cdots C_k (u,v)$, if there exists a subset of copulas $C_{k_1}\cdots C_{k_s}, $ $s\leq k\in \mathbb{N}$ such that $\psi'(\hat{C})>0 \quad \text{for}\quad \hat{C}=C_{k_1}*\cdots*C_{k_s},$ then $\psi'_{s}(C)>0$ and any Markov chain generated by $$C=a_1C_1+\cdots+a_k C_k\quad \text{for } \quad 0<a_1,\dots,a_k<1 \quad \text{is exponential}\quad \psi'-\text{mixing}.$$ 
\end{theorem}
\begin{theorem}\label{TheoL3}
For any set of copulas $C_1(u,v)\cdots C_k (u,v)$, if there exists a subset of copulas $C_{k_1}\cdots C_{k_s}, $ $s\leq k\in \mathbb{N}$ such that the density of the absolutely continuous part of $\hat{C}(u,v)$ is bounded away from $0$ $\text{for}\quad \hat{C}=C_{k_1}*\cdots*C_{k_s},$ then $\psi'_{s}(C)>0$ and any Markov chain generated by $$C=a_1C_1+\cdots+a_k C_k\quad \text{for } \quad 0<a_1,\dots,a_k<1 \quad \text{is exponential}\quad \psi'-\text{mixing}.$$ 
\end{theorem}

\subsection{All about $\psi$-mixing and $\psi^*$-mixing}
It's been shown in the literature that $\psi$-mixing implies $\psi'$-mixing, $\psi^*$-mixing and other mixing conditions. See for instance Bradley (2007). We  emphasize here that the above theorems cannot be extended to $\psi$-mixing in general by exhibiting cases when the conditions of the theorems are satisfied, but the $\psi$-mixing condition is not. A result of Bradley (1983) states the following.
\begin{lemma}\label{lemma1}
For a strictly stationary mixing sequence, either $\psi^*_n=\infty$ for all $n$ or
$\psi^*_n\rightarrow 1$ as $n\rightarrow \infty$.
\end{lemma} 
Based on this finding, if we want to show that a stationary Markov chain is $\psi^*-mixing$, it
is enough to show that it is mixing and $\psi^*_1\neq\infty$. It needs to be clear that this is not a necessary condition. In fact, there is $\psi^*$-mixing whenever we can show that for some positive integer $n$, $\psi^*_n\ne \infty$.  A remark of Longla and al. (2021) states the following.
\begin{remark} In general, for any convex convolution of two copulas (here $0 \leq a \leq 1)$,
the $\psi-mixing$ coefficient satisfies the following inequalities:\begin{eqnarray}\label{rmk}
\psi(aC_1 + (1-a)C_2)&\leq & a\psi(C_1) + (1-a) \psi(C_2);\\
\psi(aC_1 + (1-a)C_2)& \geq & a\psi(C_1) -(1-a) \psi(C_2).\end{eqnarray}
\end{remark}
A result of Longla and al. (2021) states the following.
\begin{theorem}\label{theo2}
A convex combination of copulas generates stationary $\psi-mixing$ Markov
chains if each of the copulas of the combination generates $\psi-mixing$ stationary Markov
chains.
\end{theorem} 
This Theorem as stated was not fully proved. Based on the provided proof, the correct statement should be the following.
\begin{theorem}\label{theo2}
A convex combination of copulas generates stationary $\psi-mixing$ Markov
chains if each of the copulas of the combination generates $\psi-mixing$ stationary Markov
chains with $\psi_1< 1$.
\end{theorem}

We now state the following new result.
\begin{theorem}\label{dbounded}
If a copula $C(u,v)$ is absolutely continuous and for some positive integer $n$, the density of $C^n(u,v)$ is bounded above on $[0,1]^2$, then it generates $\psi^*$-mixing stationary Markov chains. Alternatively, if for every $n$ the density of $C^n(u,v)$ is continuous and not bounded above on some subset of $[0,1]^2$, then $C(u,v)$ doesn't generate $\psi^*$-mixing or $\psi$-mixing Markov chains.
\end{theorem}

\subsubsection{Examples}
\begin{enumerate}
\item The bivariate \textit{Gaussian} copula and the Markov chains it generates. 

The Bivariate Gaussian Copula density is defined as $$c_R(u,v)=\frac{1}{\sqrt{|R|}}e^{-\frac{1}{2}(\Phi^{-1}(u)_{} \quad{}_{} \Phi^{-1}(v))(R^{-1}-\mathbb{I}){\Phi^{-1}(u)\choose \Phi^{-1}(v)}},$$
where $R$ is a bivariate variance-covariance matrix and $\mathbb{I}$  is the $2\times2$ identity matrix and $\Phi^{-1}(x)$ is the quantile function of the standard normal distribution.  The example when $R={2 \quad 1 \choose 1 \quad 1}$ is $$c_R(u,v)=e^{\Phi^{-1}(u)\Phi^{-1}(v)-.5(\Phi^{-1}(v))^2}.$$ It is clear that this density is not bounded above because for $v=.51$ and $u\to 1$, we have $c_R(u,.51)\to \infty$. By simple computations, we can show that any bivariate Gaussian copula that is not the independence copula has a density that is not bounded above. And a $*$-product of Gaussian copulas is the independence copula only when one of the two copulas is the independence copula. This is important for the following clain.
\begin{lemma}\label{lem1}
Any Copula-based Markov chain generated by a Gaussian copula that is not the product copula is not $\psi^*$-mixing.
\end{lemma}
The proof of Lemma \ref{lem1} is an application of Theorem \ref{dbounded} and the fact that the joint distribution of $(X_0,X_n)$ is the consecutive $*$-product of Gaussian copulas.

\item The Ali-Mikhail-Haq copula and the Markov chains they generate.

Copulas from the Ali-Mikhail-Haq family are defined for $\theta\in[-1,1]$ by $$C_\theta(u,v)=\frac{uv}{1-\theta(1-u)(1-v)} \quad \text{with density}\quad  c_\theta(u,v)=\frac{(1-\theta)(1-\theta(1-u)(1-v))+2\theta uv}{(1-\theta(1-u)(1-v))^3}.$$
It is easy to see that this density is continuous and satisfies $c_{\theta}(u,v)\leq \frac{1+\theta^2}{(1-\theta)^3}$ when $1>\theta>0$ or $c_{\theta}(u,v)\leq 1+\theta^2$ when $\theta\leq 0$. Therefore, the following result follows from Theorem \ref{dbounded}.
\begin{lemma}
Any copula from the Ali-Mikhail-Haq family of copulas with $\theta\ne 1$ generates $\psi^*$-mixing stationary Markov chains.
\end{lemma}
\item Copulas with densities $m_1, m_2, m_3$ and $m_4$ of Longla (2014) and the Markov chains they generate.

Because each of these copulas is bounded when the functions $g(x)$ and $h(x)$ used in their definitions are bounded, we have the following result.  
\begin{lemma}
All copulas with densities $m_1, m_2, m_3$ and $m_4$ of Longla (2014) with bounded functions $g(x)$ and $h(x)$ generate $\psi^*$-mixing Markov chains.
\end{lemma}
\end{enumerate}

\subsubsection{The Farlie-Gumbel-Morgenstern copula Family}
This family of copulas is defined by $C_{\theta}(u,v)=uv+\theta uv(1-u)(1-v)$, for $\theta\in[0,1]$. 
\begin{theorem}\label{FGM}
For any member of the Farlie-Gumbel-Morgenstern family of copula with parameter $\theta$, the joint distribution of $(X_0,X_n)$ for a stationary copula-based Markov chain generated is 
\begin{equation}
C_{\theta}^n(u,v)=uv+3\large(\frac{\theta}{3}\large)^n uv(1-u)(1-v).
\end{equation}
\end{theorem}

The density of this copula is $c^n_{\theta}(u,v)=1+3\large(\frac{\theta}{3}\large)^n(1-2u)(1-2v)$. Via simple calculations, il follows that 
\begin{equation}
0\leq 1-3\large(\frac{|\theta|}{3}\large)^n\leq c^n_{\theta}(u,v)\leq 1+3\large(\frac{|\theta|}{3}\large)^n. \label{FGM}
\end{equation}
\begin{theorem}
Any Copula-based Markov chain generated by a copula from the Farlie-Gumbel-Morgenstern family of copulas is $\psi$-mixing (for any $\theta\in[-1,1]$).
\end{theorem}
It has been established, using the first inequality of \eqref{FGM} when $n=1$ and a weaker form of Theorem \ref{TheoL3}, that any copula from this family with $|\theta|\ne 1$ generates exponential $\psi'$-mixing.  We can now show via integration that for any copula-based Markov chain $(X_1,\cdots, X_k)$ generated by $C_{\theta}(u,v)$, if $A\in\sigma(X_1)$ and $B\in\sigma(X_{n+1})$, then 
\begin{equation}
1-3\large(\frac{|\theta|}{3}\large)^n \leq \frac{P^n(A\cap B)}{P(A)P(B)}\leq 1+3\large(\frac{|\theta|}{3}\large)^n. \label{bound}
\end{equation} 
Formula \eqref{bound} implies that $\displaystyle \sup_{A,B}\frac{P^n(A\cap B)}{P(A)P(B)}\leq 1+3\large(\frac{|\theta|}{3}\large)^n<2$, for $n> 1$ and any $|\theta|\le 1$.
It follows from Theorem 3.3 of Bradley (2005) that this Markov chain is exponential $\psi$-mixing for all values of $\theta$ in the range.

\subsubsection{The Mardia and Frechet Families of Copula}
Any copula from the Mardia family is represented as $\displaystyle C_{\alpha, \beta}(u,v)=\alpha M(u,v)+\beta W(u,v)+ (1-\alpha-\beta)\Pi(u,v),$ with $0\leq \alpha, \beta, 1-\alpha-\beta \leq1$. The Frechet family of copulas is a subclass of the Mardia family with $\alpha+\beta=\theta^2$. The two families thus enjoy the same mixing properties and their analysis is theoretically identical. The density of any copula of these families is bounded away from zero on a set of Lebesgues measure 1. Therefore, the results of this paper imply that these families generate $\psi'$-mixing. Now, Consider $(X_1,X_2)$ with joint distribution $C_{\alpha,\beta}(u,v)$ and the sets $A=(0,\varepsilon)$ and $B=(1-\varepsilon, 1)$. Via simple calculations, we obtain
\begin{equation}
P(A\cap  B)=(1-\alpha-\beta)\varepsilon^2+\beta\varepsilon.
\end{equation}
Thus, 
\begin{equation}
\sup_{A,B}\frac{P(A\cap B)-P(A)P(B)}{P(A)P(B)}\ge sup_{\varepsilon}(-\alpha-\beta+\frac{\beta}{\varepsilon})=\infty. \label{nopsi}
\end{equation}
To complete the proof, we use the fact that based on the result of Longla (2014), the joint distribution of $(X_1, X_{n+1})$ is $C^n(u,v)$ - member of the Mardia family of copulas. This fact and formula \eqref{nopsi} imply that $\psi_n=\infty$ for all $n$. Therefore, this copula doesn't generate $\psi$-mixing and therefore as a result of Lemma \ref{lemma1}. Hence, the results of this work cannot be extended to $\psi$-mixing. The idea of this proof leads to the following.
\begin{theorem}\label{psi}
Let $C(u,v)$ be a copula that generates non $\psi^*$-mixing stationary Markov chains. Any convex combination of copulas containing $C(u,v)$ generates non $\psi^*$-mixing Markov chains.
\end{theorem}
Theorem \ref{psi} combined with Longla and al (2022) imply the following result.
\begin{theorem}\label{psiconvex}
A convex combination of copulas generates $\psi^*$-mixing stationary Markov chains if every copula it contains generates $\psi^*$-mixing stationary Markov chains with $\psi^*_1<1$.
\end{theorem}

\subsubsection{General case of lack of $\psi$-mixing in presence of $\psi'$-mixing}
We want here to present a large class of copulas that generate $\psi'$-mixing Markov chains, but doesn't generate $\psi^*$-mixing or $\psi$-mixing Markov chains. Based on the results of this work, we can state the following general corollary.
\begin{corollary}\label{Theo4}
Any convex combination of copulas containing the independence copula $\Pi(u,v)$ and $M(u,v)$  or $W(u,v)$ generates exponential $\psi'$-mixing, but doesn't generate $\psi$-mixing or $\psi^*$-mixing stationary Markov chains. 
\end{corollary}

\section{Some graphs of copulas and their perturbations}
In this section, we provide graphical representations of the impact of perturbations on Markov chains generated by the copulas of interest. The case is presented for some examples from the Frechet and Farlie-Gumbel-Morgenstern families of copulas. Examples are chosen for the values of the parameters that are close to independence and the extreme case of each of the families.  Two graphs of data on $(0,1)^2$ are provided as well as two graphs for the standard mornal distribution as marginal distribution of the Markov chains. To generate a Markov chain with a copula from this family, we proceed as follows.
\begin{description}
\item (a) Generate $U_1$ from $Uniform (0,1)$; 
\item (b) For $t=2,\cdots n,$ generate $W_t$ from $Uniform(0,1)$ and solve for $U_t$ the equation
$W_t= U_t+\theta (1-2U_{t-1})U_t(1-U_{t})$;
\item (c) $Y_t=G^{-1}(U_t)$, where $G(t)$ is the common marginal distribution of the variables of the stationary Markov chain.
\end{description}

Longla and al. (2021) worked on perturbation of copulas and their perturbations. For a copula $C(u,v)$, some of the studied perturbations are as follows. Assume $\alpha,\theta\in[o,1]$.

\begin{eqnarray}
\tilde{C}_{\alpha}(u,v)&=&C(u,v)+\alpha\left(\Pi(u,v)-C(u,v)\right),\label{pi}\\
\hat{C}_{\alpha}(u,v)&=&C(u,v)+\alpha\left(\text{M}(u,v)-C(u,v)\right).\label{M}
\end{eqnarray}
Formulas \eqref{pi} and \eqref{M} lead to the following.

\begin{proposition} Let $\theta\in[0,1]$, $\alpha\in[-1,1]$ and $C_{\theta}(u,v)$ be a Farlie-Gumbel-Morgenstern copula. 
\begin{eqnarray}
\tilde{C}_{\alpha,\theta}(u,v)&=&C_{\theta}(u,v)+\alpha\left(\Pi(u,v)-C_{\theta}(u,v)\right);\label{perturb11}\\
\hat{C}_{\alpha,\theta}(u,v)&=&C_{\theta}(u,v)+\alpha\left(\text{M}(u,v)-C_{\theta}(u,v)\right).\label{perturb22}
\end{eqnarray}
\begin{enumerate}
\item $\tilde{C}_{\alpha,\theta}(u,v)=C_{\theta(1-\alpha)}(u,v)$ - is a member of the Farlie-Gumbel-Morgenstern family of copulas and generates $\psi$-mixing Markov chains.
\item $\hat{C}_{\alpha,\theta}(u,v)$ is not a member of the Farlie-Gumbel-Morgenstern family of copulas and does not generates $\psi-mixing$ Markov chains.
\end{enumerate}
\end{proposition}

\begin{figure}[ht!]
\includegraphics[width=17cm, height=8.5cm]{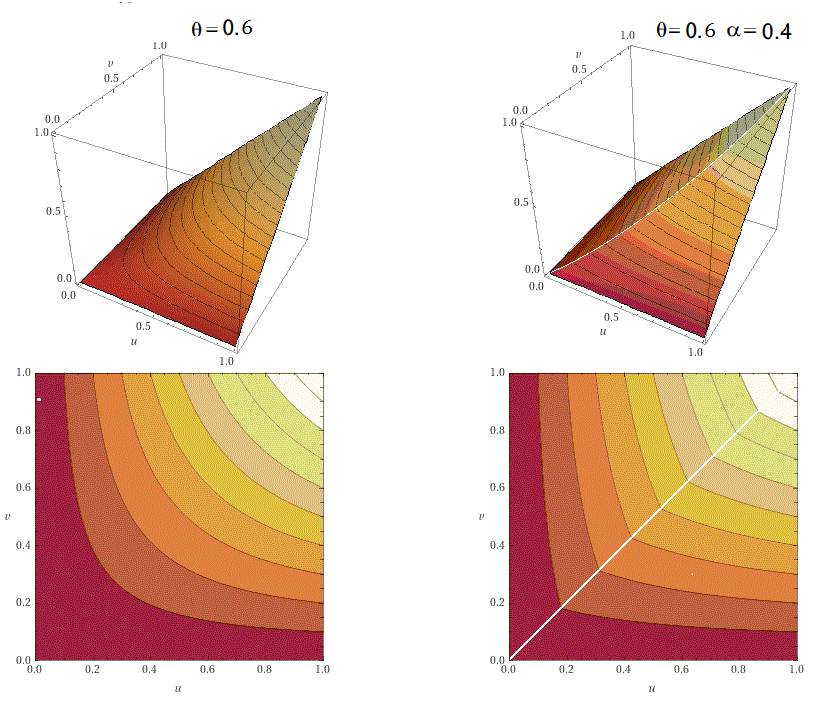}
\caption{Farlie-Gumbel-Morgenstern copula and level curves}\label{fig1}
\end{figure}

\begin{figure}[!ht]
\begin{center}
\includegraphics[width=18cm, height=10cm]{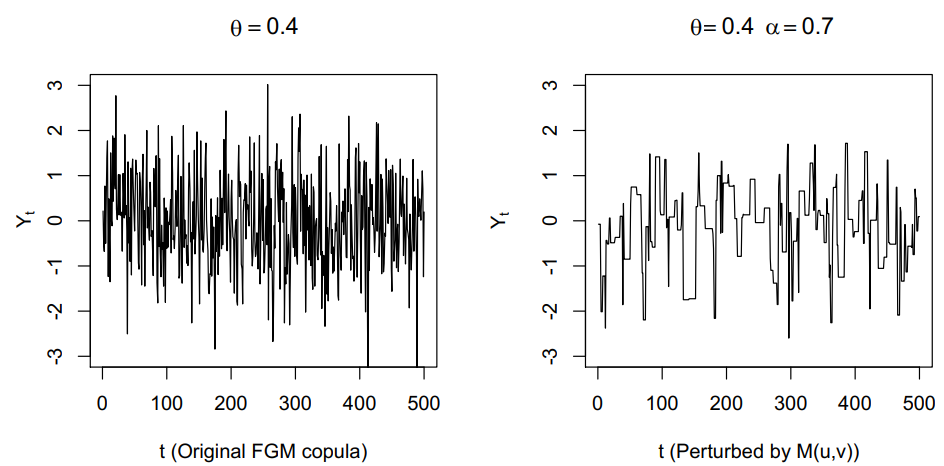}
\caption{\it Data from the Farlie-Gumbel-Morgenstern copula and its perturbations.}\label{fig4}
\end{center}
\end{figure}

On Fifure \ref{fig1} we have a 3-dimensional graph of the Farlie-Gumbel-Morgenstern copula with parameter $\theta=.6$ and its level curves on the left and the corresponding graphs for the perturbation with parameter $\alpha=.4$ on the right. Figure \ref{fig4} represents a simulated Markov chain form the Farlie-Gumbel-Morgenstern copula with $\theta=.4$ and the one generated by its perturbation with parameter $\alpha=.7$. Here, the marginal distribution of the Markov chain is standard normal. We can see on the graphs that the mixing structure is not the same when the copula is perturbed by $M(u,v)$. This supports the theoretical results. 

The Mardia family of copulas is defined by 
\begin{equation}\label{Mardia}
C_{a,b}(u,v)=aM(u,v)+bW(u,v)+(1-a-b)\Pi(u,v)
\end{equation}
and the Frechet copulas are a subfamily with $a=\dfrac{\theta^2(1+\theta)}{2}$, $ b=\dfrac{\theta^2(1-\theta)}{2}$ and $ |\theta|\leq1$.
Unlike Farlie-Gumbel-Morgenstern copulas, these copulas are not absolutely continuous.  To generate an  observation $(U,V)$ from $C_{\theta}(u,v)$, one needs to generate independent observations $(U,V_1, V_2)$ from the uniform distribution on $(0,1)$. Then, do the following:  

$$V=\left\{\begin{array}{lcl} V_2 & \text{if}& V_1< 1-\theta^2,\\
 U & \text{if}& 1-\theta^2<V_1<1-\theta^2+\theta^2(1+\theta)/2,\\
1-U &\text{if}& V_1>1-\theta^2+\theta^2(1+\theta)/2.
\end{array}
\right.$$

\begin{figure}[ht!]\begin{center}\includegraphics[width=17cm, height=7.5cm]{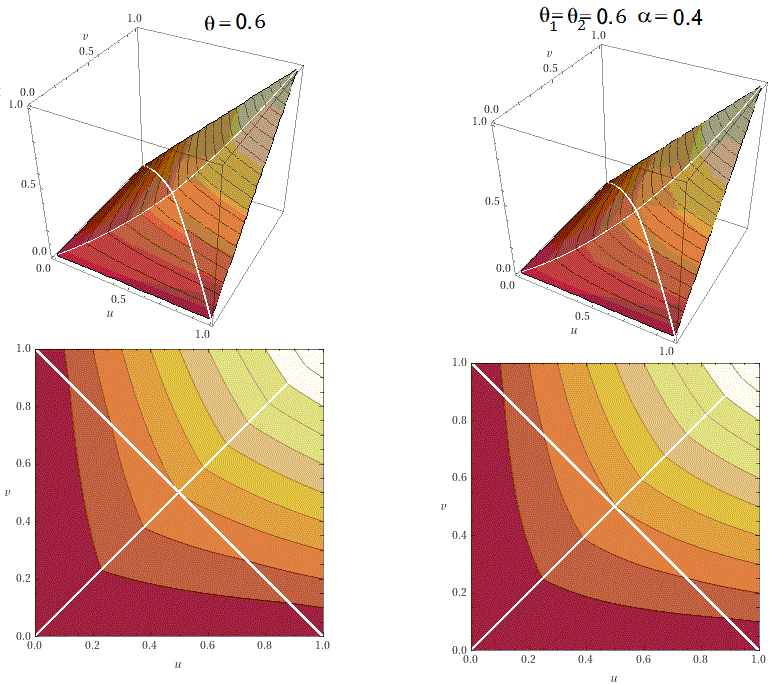}\caption{\it Frechet copula represenation and level curves}\label{frechetrep}\end{center}\end{figure}
\begin{figure}[ht!]
\begin{center}
\includegraphics[width=18cm, height=7.5cm]{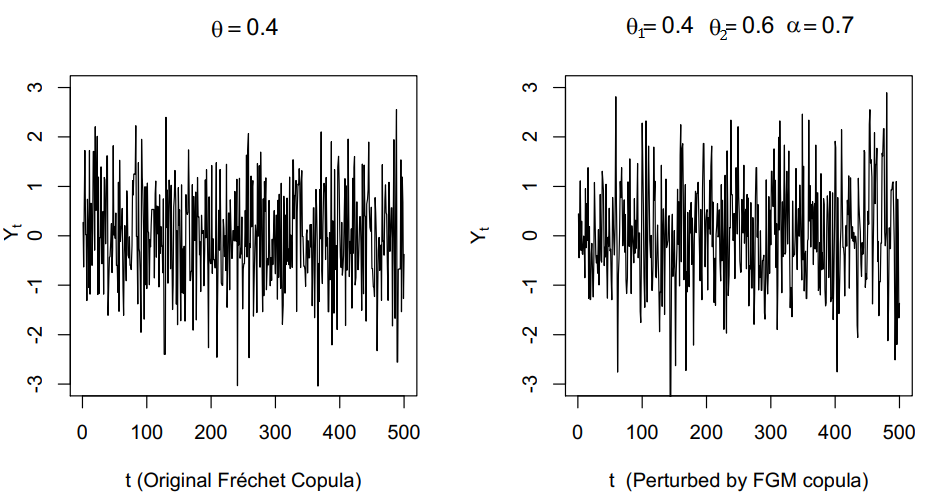}
\caption{\it Markov chain generated by Frechet copulas and its perturbations.}\label{fig3}
\end{center}
\end{figure}

 
Figure \ref{frechetrep} gives a representation of the Frechet copula for $\theta=.6$ and its perturbation with $\alpha=.4$, together with level curves. Figure \ref{fig3} represents a Markov chain with 500 observations simulated from the Frechet copula with $\theta=.6$ and its perturbation with parameter $\alpha=.7$. 
Perturbations of the Frechet copula will have the form:
\begin{eqnarray}\label{frechpert}
\tilde{C}_{\theta,\alpha}(u,v)=C_\theta(u,v)+\alpha(\Pi(u,v)-C_\theta(u,v));\label{frechpert1}\\
\hat{C}_{\theta,\alpha}(u,v)=C_\theta(u,v)+\alpha(M(u,v)-C_\theta(u,v)).\label{frechpert2}
\end{eqnarray}

It is good to notice that these perturbations are not Frechet copulas, but remain in the class of Mardia copulas.  Figure \ref{fig3} represents a Markov chain generated by a Frechet copula and the ones generated by its perturbations using the standard normal distribution for marginal distributions.

\section{Simulation study}
This simulation study shows the importance of this topic. We simulate a dependent data set that exhibits $\psi$-mixing or $\psi'$-mixing and show how the mixing structure influences the statistical study. Based on the fact that the considered mixing coefficient converges exponentially to $0$, we can bound the variance o partial sums and obtain the condition of the central limit theorem and confidence interval of Longla and Peligrad (2020). Thanks to this central limit theorem, we construct confidence intervals without having to estimate the limiting  variance of the central limit theorem of Kipnis and Varadhan (1986) that holds here because the Markov chains are reversible and $n var(\bar{Y})\to \sigma<\infty$. The standard central limit theorem is useless in this case because the limiting variance is not necessarily that of $Y$. Let us recall here the formulations of Longla and Peligrad (2020). They have proposed a new robust confidence interval for the mean based on a sample of dependent observations with a mild condition on the variance of partial sums. This confidence interval needs a random sample $(X_i, 1\leq i\leq n)$, generated independently of $(Y_i, 1\leq i\leq n)$ and following the standard normal distribution, the Gaussian Kernel and the optimal bandwidths 
$$h_n=\left[\dfrac{\bar{y^2_n}}{n\sqrt{2}\bar{y}^2_n}\right]^{1/5}.$$

Let us check the conditions required for use of this proposed estimator of the mean and the confidence interval.  They are as follows: 
\begin{enumerate}
\item  $(Y_i)_{i\in\mathbb{Z}}$ is an ergodic sequence;
\item  $(Y_i)_{i\in\mathbb{Z}}$ have finite second moments;
\item  $nh_n var(\bar{Y}_n)\rightarrow 0$ as $n\rightarrow \infty$.
\end{enumerate}
For the sake of clarity, we will use $C^{FGM}_\theta(u,v)$ to denote the Farlie-Gumbel-Morgenstern copula with parameter $\theta$.

\textbf{Verification of the conditions}

\begin{enumerate}
\item Ergodicity
\begin{enumerate}
\item It has been shown by Theorem 2.3 and Example 2.4 of Longla (2014) that the copula $C_{\theta}^{FGM}(u,v)$ generates geometrically ergodic Markov chains.
\item From this current work, we can deduce that the perturbed $\hat{C}^{FGM}_{\theta,\alpha}(u,v)$ generates $\psi'-mixing$ Markov chains. In fact, this copula is a linear combination of two copulas such that one is $\psi'-mixing$. In addition, (see Bradley (2005) and Longla and Peligrad (2012)) $\psi'-mixing$ implies $\phi-mixing$ and $\phi-mixing$ implies geometric ergodicity for reversible Markov chains. So the Markov chain generated by $\hat{C}^{FGM}_{\theta,\alpha}(u,v)$  is geometrically ergodic.
\item According to Theorem 2.16 and Remark 2.17 of Longla (2014), the Frechet copula $C_\theta(u,v)$ generates geometrically ergodic Markov chains.
\item  The  perturbed Frechet copula $\hat{C}_{(\theta_1,\theta_2,\alpha)}(u,v)$ is a linear combination of copulas $C_{\theta_1}(u,v)$ and $C^{FGM}_{\theta_2}(u,v)$. These two copulas are  symmetric  and each one generates geometrically ergodic sequences as said above. Then, according to Theorem 5 of Longla and Peligrad (2012), this copula generates geometrically ergodic Markov chains.
\end{enumerate}
\item  The marginal distribution used in this work is normal with mean 30 and variance 1. Therefore, it has second moments.
\item The condition on the variance ($nh_n var(\bar{Y})\to 0$) is checked in the appropriate section below.
\end{enumerate}

For data simulation, we set $Y_i\sim N(30,1)$ for all copulas and the perturbation parameter $\alpha=0.4$ in all cases. For Farlie-Gumbel-Morgenstern and Frechet copulas we set $\theta=0.6$. For the Frechet perturbed copula, $\theta_1=\theta_2=0.6$. 
For $1\leq i\leq n$, $X_i\sim N(0,1)$ is a sequence of independent random variables that is independent of the Markov chain $(Y_i, 1\leq i\leq n)$.

According to above considerations, the estimator of $\mu_Y$ is $\tilde{r}_n=\dfrac{1}{nh_n}\sum\limits_{i=1}^nY_i \exp\left(-0.5(\dfrac{X_i}{h_n})^2\right)$ and the confidence interval is $\left(\tilde{r}_n\sqrt{1+h_n^2}-z_{\alpha/2}\left(\dfrac{\bar{Y_n^2}}{nh_n\sqrt{2}}\right)^{1/2};\tilde{r}_n\sqrt{1+h_n^2}+z_{\alpha/2}\left(\dfrac{\bar{Y_n^2}}{nh_n\sqrt{2}}\right)^{1/2}\right)$.\\

The following table is the result of the simulation study for the Markov chains generated by the two considered copulas and their perturbations.
\begin{center}
\begin{tabular}{|c|c|c|c|c|c|}
\hline
\text{Copula } &size&n=100&n=5000&n=10000&n=20000\\
\hline
\multirow{2}{*}{$C^{FGM}_{\theta}$}&\text{Estimator of } $\mu_Y$&23.25&28.41&29.85&29.54\\
\cline{2-6}
&\text{Confidence interval} &(16.72, 32.88)&(27.12, 30.51)&(28.89, 31.46)&(28.81, 30.76)\\
\hline
\multirow{2}{*}{$\hat{C}^{FGM}_{\theta,\alpha}$}&\text{Estimator of } $\mu_Y$&23.70&28.39&29.80&29.52\\
\cline{2-6}
&\text{Confidence interval} &(17.08, 33.50)&(27.10, 30.50)&(28.84, 31.42)&(28.79, 30.74)\\
\hline
\multirow{2}{*}{$C_{\theta}$}&\text{Estimator of } $\mu_Y$&31&29.40&30.30&30.29\\
\cline{2-6}
&\text{Confidence interval} &(24.97, 41,16)&(28.13, 31.52)&(29.34,  31.91)&(29.56, 31.51)\\
\hline
\multirow{2}{*}{$\hat{C}_{(\theta_1,\theta_2,\alpha)}$}&\text{Estimator of } $\mu_Y$&31.15&29.39&30.29&30.23\\
\cline{2-6}
&\text{Confidence interval} &(25.08, 41.37)&(28.11, 31.51)&(29.33, 31.90)&(29.51, 31.46)\\
\hline
\end{tabular}
\end{center}

\section{Conclusion and remarks} The graphs and simulations presented in this paper have been obtained using $R$. We have provided some insights on $\psi^*$-mixing, $\psi'$-mixing and $\psi$-mixing. Though we have presented extensive examples and results for $\psi'$-mixing and $\psi^*$-mxing, we have not been able to answer the question on convex combinations of $\psi$-mixing. The following question remains open: Does a convex combination of $\psi$-mixing generating copulas generate $\psi$-mixing? A positive answer to this question has been presented for the case when each of the copulas satistfy $\psi_1<1$. It would also be interesting to find a general condition on the copula for $\psi$-mixing  like the one presented for $\psi^*$-mixing.
 
\section{Appendix of proofs}
\subsection{Proof of Theorem \ref{TheoL}}

Recall that the function $c(x,y)$ defined on $I^2$ is said to be bounded away from zero on a set of Lebesgue measure 1 iff
$\exists m>0, m\in \mathbb{R}, \exists Q\subset I^2: \lambda(Q)=1, \forall (x,y)\in Q$,  $c(x,y)\geq m.$

According to Theorem \ref{Brad}  a strictly stationary Markov chain $(X_k, ~ k\in\mathbb{N})$ is $\psi'$-mixing if :  $
\text{ for some } n\in\mathbb{N}, ~ \psi_n'(C)\neq 0.$
Let $A\subset I,~ B\subset I$, by easy calculation we obtain:
\begin{eqnarray}
\dfrac{\int_A\int_Bc_1(x,y)dxdy}{\lambda(A)\lambda(B)}&=&\dfrac{\int_A\int_Bc(x,y)dxdy}{\lambda(A)\lambda(B)}\nonumber\\
&\geq&\dfrac{\int_A\int_B(\varepsilon_1(x)+\varepsilon_2(y))dxdy}{\lambda(A)\lambda(B)}\nonumber\\
&=&\dfrac{\int_A\int_B \varepsilon_1(x)dxdy}{\lambda(A)\lambda(B)}+\dfrac{\int_A\int_B\varepsilon_2(y)dxdy}{\lambda(A)\lambda(B)}\nonumber\\
&=&\dfrac{\int_A\varepsilon_1(x)dx \int_B dy}{\lambda(A)\lambda(B)}+\dfrac{\int_B\varepsilon_2(y)dy\int_Adx}{\lambda(A)\lambda(B)}\nonumber\\
&=&\dfrac{ \lambda(B)\int_A\varepsilon_1(x)dx}{\lambda(A)\lambda(B)}+\dfrac{\lambda(A)\int_B\varepsilon_2(y)dy}{\lambda(A)\lambda(B)}
\end{eqnarray}
So for all $A\subset I,~ B\subset I$, the following inequality holds:
\begin{equation}\label{lab}
\dfrac{\int_A\int_Bc_1(x,y)dxdy}{\lambda(A)\lambda(B)}\geq \dfrac{ \int_A\varepsilon_1(x)dx}{\lambda(A)}+\dfrac{\int_B\varepsilon_2(y)dy}{\lambda(B)}
\end{equation}
We also have:
\begin{center}
$ \dfrac{ \int_A\varepsilon_1(x)dx}{\lambda(A)}\geq \underset{\underset{\lambda(A)>0}{A\subset I}}{\inf}\dfrac{ \int_A\varepsilon_1(x)dx}{\lambda(A)} \text{ and } \dfrac{\int_B\varepsilon_2(y)dy}{\lambda(B)}\geq \underset{\underset{\lambda(B)>0}{A\subset I}}{\inf}\dfrac{\int_B\varepsilon_2(y)dy}{\lambda(B)}.$
\end{center}
Thus,  for all $A\subset I,~ B\subset I$ :\begin{center}
$\dfrac{\int_A\int_Bc_1(x,y)dxdy}{\lambda(A)\lambda(B)}\geq  \underset{\underset{\lambda(A)>0}{A\subset I}}{\inf}\dfrac{ \int_A\varepsilon_1d\lambda}{\lambda(A)} + \underset{\underset{\lambda(B)>0}{A\subset I}}{\inf}\dfrac{\int_B\varepsilon_2d\lambda}{\lambda(B)}.$
\end{center}
Which means:
\begin{center}
$\underset{\underset{ \lambda(A)\lambda(B)>0}{ A\subset I,~ B\subset I}}{\inf}\dfrac{\int_A\int_Bc_1(x,y)dxdy}{\lambda(A)\lambda(B)}\geq M+N.$
\end{center} Where $M= \underset{\underset{\lambda(A)>0}{A\subset I}}{\inf}\dfrac{ \int_A\varepsilon_1d\lambda}{\lambda(A)}$ and $N= \underset{\underset{\lambda(B)>0}{A\subset I}}{\inf}\dfrac{\int_B\varepsilon_2d\lambda}{\lambda(B)}$.
Hence, 
$
\psi'_1(C)\geq M+N.$
According to the theorem assumptions, $M >0$ or $ N>0$.  So,
$
\psi'_1(C)\geq M+N>0.
$
We can conclude $(X_k, k\in\mathbb{N})$ is $\psi'$-mixing.

\subsection{Theorem \ref{TheoL2} and Theorem \ref{TheoL3}}
To prove these theorems, we will use the following proposition from Longla and al (2022)

\begin{proposition}\label{prop1}
For a convex combination of copulas  $\displaystyle C(x,y)=\sum_{i=1}^ka_i C_i(x,y),$ where $0<a_1,...,a_k<1$ and $\displaystyle\sum_{i=1}^k a_i=1$,
 the following formula holds. For any $s\in\mathbb{N}$
 \begin{equation}\label{e1}
 C^s(x,y)=\sum_{j=1}^{k^s}b_j\times ~_1C_j\ast...\ast~_sC_j(x,y),
 \end{equation}
 where $\sum_{j=1}^{k^s}b_j=1,~ 0<b_1,..., b_{k_s}<1 $, and each of the copulas $~_iC_j(x,y)=C_{j_i}(x,y)$ for some $j_i\in\{1,..., k\}$ and the sum is over all possible products of $s$ copulas selected from the original $k$ copulas with replacement.
\end{proposition}
The notation $~_iC_j$ indicates that the copula $C_{j_i}$ was selected in the given $j^{th}$ element of $B= \{C_1,...,C_k\}^n$.
  
\textbf{(1) ~} Suppose that  there exists a subset of copulas $C_{k_1},...,C_{k_s}~, s\leq k\in\mathbb{N}$ such that $\psi'(\hat{C})>0$ for $\hat{C}=C_{k_1}\ast ...\ast C_{k_s}$. Equation  (\ref{e1}) can be written as follows:
\begin{equation}\label{e2}
C^s(x,y)=b_i \hat{C}(x,y)+\sum_{\underset{j\neq i}{j=1}}^{k^s}b_j \hat{C}_j(x,y),
~~~\text{where}\end{equation} 
 $\hat{C}(x,y)=C_{i_1}\ast...\ast C_{i_s}(x,y)$~~ and  $\hat{C}_j(x,y)=C_{j_1}\ast...\ast C_{j_s}(x,y)$ \\
 
 Let $(X_k, k\in\mathbb{N})$ be a copula-based Markov chain generated by the copula $C(x,y)$; $(\hat{X}^j_k, k\in \mathbb{N})$ a Markov chain generated by copula $\hat{C}_j$ for $1\leq j \leq k^s$, $\hat{C}_i=\hat{C}$. For $A\in\sigma(X_0)$ and $B\in\sigma(X_s)$, equation (\ref{e2}) yields
 \begin{eqnarray}\label{e3}
P^s(A\cap B)&=&b_i\hat{P}(A\cap B)+ \sum_{\underset{j\neq i}{j=1}}^{k^s}b_j \hat{P}_j(A\cap B)\ge b_i \hat{P}(A\cap B),
 \end{eqnarray}
 where $P^s(A\cap B)= P(X_1\in A, X_{s+1}\in B)$; $\hat{P}_j(A\cap B)=P(X^j_1\in A, X^j_{s+1}\in B)$ and $\hat{P}(A\cap B)=P(X^i_1\in A, X^i_{s+1}\in B)$. Thus, 
 \begin{center}
 $\psi'_s(C)=\underset{A\subset I, B\subset I, P(A)(B)>0~}{\inf}\dfrac{P^s(A\cap B)}{P(A)P(B)}\ge b_i\psi'(\hat{C})$.
 \end{center}
 By our assumptions,  $\psi'(\hat{C})>0$. The conclusion follows from Theorem \ref{Brad}.
 
\textbf{ (2)}~ Suppose there exists a subset of copulas $C_{k_1},...,C_{k_s}~, s\leq k\in\mathbb{N}$ such that the density of the absolutely continuous part of the copula  $\hat{C}=C_{k_1}\ast ...\ast C_{k_s}$ is bounded away from zero. From equation  (\ref{e2}) we have:
\begin{equation}\label{f1}
c^s(x,y)\geq b_i \hat{c}(x,y).\end{equation} 
Moreover, the density of the absolutely continuous part of $\hat{C}(u,v)$ is bounded away from zero. Thus, there exists $c>0$: $\forall (x,y)\in[0,1]^2$, ~$\hat{c}(x,y)\geq c$ almost surely. Hence, from (\ref{f1}), we have
$c^s(x,y)\geq b_i c$.  Now, if $(X_k, k\in\mathbb{N})$ is a copula-based Markov chain generated by the copula $C(x,y)$ and an absolutely continuous distribution, then for $A\in\sigma(X_1)$ and $B\in\sigma(X_{s+1})$, we have
 \begin{eqnarray}\label{f4}
P^s(A\cap B) \geq  b_i c P(A)\times P(B)
\quad \text{ and } \quad \dfrac{P^s(A\cap B)}{P(A)\times P(B)} \geq b_i c ,
 \end{eqnarray}
 where $P^s(A\cap B)= P(X_1\in A, X_{s+1}\in B)$. It follows from equation (\ref{f4}) that 
 $$\displaystyle \psi'_s(C)=\underset{P(A)(B)>0~}{\inf}\dfrac{P^s(A\cap B)}{P(A)P(B)}\geq b_i c >0.$$
 This concludes the proof of Theorem \ref{TheoL3}.
\subsection{Proof of Theorem \ref{FGM}}
The following decomposition if true for Farle-Gumbel-Morgenstern copulas with $\lambda=1-\theta$: $$C_{\theta}(u,v)=(1-\lambda) (uv+uv(1-u)(1-v))+\lambda uv.$$
Given that $C(u,v)=(uv+uv(1-u)(1-v))$ is a copula, we can apply Theorem \ref{copL} to obtain
$$C^n_\theta(u,v)=(1-\lambda)^nC^n(u,v)+(1-(1-\lambda)^n)uv.$$
It remains to show that $C^n(u,v)=uv+3\large(\frac{1}{3}\large)^n uv(1-u)(1-v)$ by mathematical induction. It is clear that the formula is correct for $n=1$. Assume that for $n=k$, we have 
$$C^k(u,v)=uv+3\large(\frac{1}{3}\large)^k uv(1-u)(1-v).$$
Using the fold product, we obtain $$C^{k+1}=C^{k}*C(u,v)=\int_0^1C^k_{,2}(u,t)C_{,1}(t,v)dt.$$
$$C^k_{,2}(u,t)=u+3\large(\frac{1}{3}\large)^k u(1-u)(1-2t) \quad \text{and} \quad C_{,1}(t,v)=v+v(1-v)(1-2t).$$
Plugging these functions into the integral and computing yields the needed results. The proof ends by replacing $C^n(u,v)$ by its value and using $\lambda=1-\theta$.
\subsection{Proof of Theorem \ref{dbounded}}
Assume that the copula $C(u,v)$ is such that for all $(u,v)\in [0,1]^2$, $c^n(u,v)\leq K$, where $K$ is a constant and $c^n$ is the density of the copula $C^n(u,v)$. Let $A\in \sigma(X_0)$ and $B\in\sigma(X_n)$, where $(X_0,X_n)$ has copula $C^n(u,v)$. Assume that the stationary distribution of the Markov chain has distribution $F(x)$. 
Using Sklar's Theorem (see Sklar (1959)), we have $$P^n(A\cap B)=P(X_0\in A , X_n\in B)=\int_A\int_{B}c^n(F(x), F(y))dF(y)dF(x).$$
Therefore, $P^n(A \cap B)\leq KP(A)P(B)$. This implies $\psi_1(C)\leq K\ne \infty$ and by Lemma \ref{lemma1}, $C(u,v)$ generates stationary $\psi$-mixing Markov chains. Now, if we assume that There exists a set of non-zero measure $\Omega\subset [0,1]^2$ such that $A\times B\subset \Omega$, $A\in\sigma(X_0)$, $B\in\sigma(X_1)$ and the density of $C^n$ is not bounded above on $\Omega$, but bounded below by a any given non-zero real number $M$. This construction is possible due to continuity of the density of $C(u,v)$. It follows that for any constant $M$, $$ P^n(\Omega)\ge P^n(A\times B)\geq M \int_{A\times B}d\Pi(x,y)=MP(A)P(B).$$
It is obvious here that As $M$ grows, the size of $P(A)P(B)$ reduces as their product has to be at most $1$. From here, we obtain
$$\frac{P(A\cap B)}{P(A)P(B)}\ge M. \quad \text{This leads to }\quad \psi^*_n(C)>M.$$
Because this is true for every $M$ and every $n$, we can conclude that $\psi^*_n(C)=\infty$ and $\psi_n(C)=\infty$ for all $n$. Thus, the generated Markov chain is not $\psi^*$-mixing and not $\psi$-mixing.
\subsection{Proof of Theorem \ref{psi} and Theorem \ref{psiconvex}}
Without loss of generality the proof can be done for a convex combination of two copulas, one of which is $C(u,v)$ and doesn't generate $\psi^*$-mixing Markov chains. This is true because any convex combination of copulas can be written as a convex combination of two copulas. Now, assume that $$C_2(u,v)=\alpha C(u,v)+(1-\alpha)C_1(u,v).$$ By Lemma \ref{lemma1} , $\psi^*_n(C)=\infty$ for all $n\in\mathbb{N}$. We need to show that $\psi^*_n(C_2)=\infty$ for all $n\in\mathbb{N}$. By Longla and al (2021), there exist $b_{in}, C_{1in}(u,v)$, such that $b_{in}>0$, $\alpha^n+\sum_{i=2}^{2^n}b_{in}=1$ and
$$C_2^{n}(u,v)=\sum_{i=2}^{2^n} b_{in} C_{1in}(u,v)+\alpha^n C^n(u,v).$$
Therefore, The probability distribution $P^{n}_2$ of  $(X_1, X_{n+1})$ from the Markov chain generated by $C_2(u,v)$ and the probability distributions $P^n_{1i}$ of $(\tilde{X}_{i1},\tilde{X}_{in+1})$ for the Markov chains generated by the copulas $C_{1in}(u,v)$ satisfy the following relationship for every $A\in\sigma(X_0)$ and $B\in\sigma(X_{n+1})$:
$$P_2^{n}(A\cap B)=\sum_{i=2}^{2^n} b_{in} P^n_{1i}(A\cap B)+\alpha^n P^n(A\cap B).$$
Therefore, $P_2^n(A\cap B)\ge \alpha^nP^n(A\cap B)$. Given that $\psi^*_n(C)=\infty$ for all $n$, it follows that $\sup_{A,B}\frac{P^n(A\cap B)}{P(A)P(B)}=\infty$, leading to $$\sup_{A,B}\frac{P^n_2(A\cap B)}{P(A)P(B)}=\infty \quad \text{and} \quad \psi^*_n(C_2)=\infty \quad \text{for all}\quad n\in\mathbb{N}.$$
This concludes the proof of Theorem \ref{psi}. 
Now, to prove Theorem \ref{psiconvex}, as for the previous case, it is enough to consider a convex combination of two copulas. Assume that $C_1(u,v)$ and $C_2(u,v)$ generate each $\psi^*$-mixing (or $\psi$-mixing) stationary copula-based Markov chains qith $\psi^*_1<1$ (or $\psi_1<1$) respectively. Define $C(u,v)=\alpha C_1(u,v)+(1-\alpha)C_2(u,v)$. Once more, we will use Lemma \ref{lemma1}.  $\psi_1(C)\leq \alpha \psi_1(C_1)+(1-\alpha) \psi_1(C_2)<(\alpha)+(1-\alpha)=1$. The same argument works for $\psi_1^*(C)$.

\subsection{Checking the condition $nh_n var(\bar{Y})\to 0$.}
Given the the Markov chains that we consider here are reversible and ergodic (see Haggstrom and Rosenthal (2007), Kipnis and Varadhan (1986)), $$n var(\bar{Y})\quad \text{behaves as}\quad var(Y_0)+2\sum_{k=1}^\infty cov(Y_0, Y_k).$$
Moreover, if the series converges, then the central limit theorem holds with variance equal to its sum.
On the other side, Markov chains generated by Farlie-Gumbel-Morgenstern copulas, Frechet copulas and their considered perturbations are exponential $\psi'$-mixing. This implies that they are all exponential $\rho$-mixing. Exponential $\rho$-mixing implies convergence of the considered series. Therefore $n var(\bar{Y})\to C$, leading to $nh_n var(\bar{Y})$.

\end{document}